\documentclass[10pt]{article}

\usepackage{graphics,epsfig}
\usepackage{amsmath}

\usepackage{amsmath}
\usepackage{amsfonts}

\input epsf
\usepackage{makeidx}
\makeindex
\usepackage{amsfonts}
\def\ds{\displaystyle}

\def \b0 {{\bf 0}}

\def \ff12{\frac {1}{2}}

\newtheorem{theorem}{Theorem}[section]
\newtheorem{lemma}[theorem]{Lemma}
\newtheorem{corollary}[theorem]{Corollary}

\newtheorem{re}[theorem]{Remark}
  \newtheorem{prop}[theorem]{Proposition}

\newenvironment{proof}{\noindent {\em Proof:  }}{\hspace*{1cm}
    \hspace*{\fill}$\rule{1.2ex}{1.4ex}$\vspace{.15cm}}

\newtheorem{proposition}[theorem]{Proposition}

\begin{document}

\centerline{\Large \bf SubRiemannian geometry on the sphere
$\mathbb{S}^3$}

\medskip

\centerline {\bf Ovidiu Calin, Der-Chen Chang and Irina
Markina\footnote{The first author partially supported by a
research grant at Eastern Michigan University. The second author
partially supported by a competitive research grant at Georgetown
University. The third author is partially supported by a research
grant at University of Bergen.}}
\medskip

\noindent {\small {\bf Abstract:}  We discuss the subRiemannian
geometry induced by two noncommutative vector fields which are left invariant
on the Lie group $\mathbb{S}^3$.}
{\small {\bf MS Classification} (2000): }{\small 53C17, 53C22, 35H20}\\

\section{Introduction}
The study of step 2 subRiemannian manifolds has the Heisenberg
group as a prototype. This is a noncommutative Lie group with the
base manifold $\mathbb{R}^3$ and endowed with a nonintegrable
distribution spanned by two of the noncommutative left invariant
vector fields. This structure enjoys also the property of being a
contact structure or a CR-manifold. The study of the subRiemannian
geodesics on the Heisenberg group started with the work of Gaveau
[9]. One trend in the literature is to use the geometry of the
Heisenberg group to describe the Heisenberg Laplacian and its heat
kernel, see Beals, Gaveau, Greiner [1,2,3,4]. Later, this
structure led to generalizations of the Heisenberg group as can be
seen in Calin, Chang, Greiner [5,6] and Chang, Markina [7]. For
more fundamental issues on subRiemannian geometry, see Strichartz
[10].

\smallskip
One question is, if there is another 3-dimensional noncommutative
Lie group besides the Heisenberg group. The answer is positive and
this paper deals with it. The sphere $\mathbb{S}^3$ and the group
$SU(2)$ have these properties. More important, this is the first
compact subRiemannian manifold has been under consideration.

\smallskip
The present paper starts with an introduction to quaternions and
then defines the 3-dimmensional sphere as the set of quaternions
of length one. The quaternion group induces on $\mathbb{S}^3$ a
structure of noncommutative Lie group. This group is compact and
the results obtained in this case are very different than those
obtained in the case of the Heisenberg group, which is a
noncompact Lie group. Like in the Heisenberg group case, we
introduce a nonintegrable distribution on the sphere and a metric
on it using two of the noncommutative left invariant vector
fields. This way $\mathbb{S}^3$ becomes a subRiemannian manifold.
It is known that the group $SU(2) $ is isomorphic with the sphere
$\mathbb{S}^3$ and represents an example of subRiemannian manifold
where the elements are matrices. The main issue here is to study
the connectivity by horizontal curves and its geodesics on this
manifold. In this paper, we are using Lagrangian method to study
the connectivity theorem on ${\mathbb S}^3$ by horizontal curves
with minimal arc-length. We show that for any two points in
${\mathbb S}^3$, there exists such a geodesic joining these two
points.

\smallskip
The final version of this paper was in part written while the
first two authors visited the National Center for Theoretical
Sciences and National Tsing Hua University during June-July of
2006. They would like to express their profound gratitude to
Professors Jin Yu and Shu-Cheng Chang for their invitation and for
the warm hospitality extended to them during their stay in Taiwan.

\section{The Quaternions Group $\mathbf{Q}$}

William R. Hamilton introduced quaternions in 1873. They are
obtained as a non-commutative extension of the complex numbers.
Let $ \mathbf{Q}$ be the set of matrices $\ds
\begin{pmatrix}
  \alpha   &  \beta   \\
  -\overline{\beta}  & \overline{\alpha }
\end{pmatrix} $
with complex number entries $ \alpha , \beta \in \mathbb{C}$. The
set $\mathbf {Q}$ becomes a group with the usual matrix
multiplication law, called the {\it quaternions group}. The
identity element is $\ds \mathrm{I} =\begin{pmatrix}
  1  &  0  \\
  0 & 1
\end{pmatrix} $ and the inverse of $q = \begin{pmatrix}
  \alpha   &  \beta   \\
  -\overline{\beta}  & \overline{\alpha }
\end{pmatrix}$ is
$\ds  q^{-1} = \frac{1}{\Delta} \begin{pmatrix}
  \overline{\alpha }  &  -\beta   \\
  \overline{\beta}  & {\alpha }
\end{pmatrix}$, where $\Delta = |\alpha|^2 + |\beta|^2.$
Let
\begin{equation}\label{eq:matrixun}
\mathbf{i} = \begin{pmatrix}
  i & 0 \\
  0 & -i
\end{pmatrix}, \qquad    \mathbf{j} = \begin{pmatrix}
  0 & 1\\
  -1 & 0
\end{pmatrix}, \qquad \mathbf{k} = \begin{pmatrix}
  0 & i \\
  i & 0
\end{pmatrix}. \end{equation}
Then
\begin{eqnarray*}
\mathbf{i}^2 = \mathbf{j}^2 = \mathbf{k}^2 = -\mathrm{I}, \qquad
\mathbf{i}\mathbf{j} =
 - \mathbf{j}\mathbf{i} = \mathbf{k},\\
\mathbf{j}\mathbf{k} = - \mathbf{k}\mathbf{j} = \mathbf{i} ,
\qquad \mathbf{k}\mathbf{i} = -\mathbf{i}\mathbf{k} = \mathbf{j}.
\end{eqnarray*}
Any real number $a \in \mathbb{R}$ can be identified with $\ds a
\mathbf{i}$ by the one to one homomorphism $\varphi : \mathbb{R}
\rightarrow \mathbf {Q}$, $\varphi (a) = \begin{pmatrix}
  a & 0 \\
  0 & -a
\end{pmatrix}$.
If $\alpha = a_0 + a_1 \mathbf{i}, \beta = b_0 + b_1 \mathbf{i}
\in \mathbb{C} $, then an elementary computation yields
$$ q= \begin{pmatrix}
 \alpha  & \beta \\
  -\overline{\beta}  & \overline{\alpha }
\end{pmatrix} = a_0 \mathrm{I} + a_1  \mathbf{i} + b_0  \mathbf{j} + b_1  \mathbf{k}.$$
The set $\{ \mathrm{I}, \mathbf{i}, \mathbf{j}, \mathbf{k} \}$ is
a basis of $\mathbf {Q}$. Let $ q = a \mathrm{I} + b \mathbf{i} +
c \mathbf{j }+ d \mathbf{k} $ and $ q' = a' \mathrm{I} + b'
\mathbf{i} + c' \mathbf{j }+ d' \mathbf{k} $ be two quaternions.

A straightforward computation shows
\begin{eqnarray*}
qq' &=& (aa' - bb' -cc' - dd') \mathrm{I} +  (ab' + ba' + cd' -
dc')\mathbf{ i }\\
&&+(ac'+ca'+db'-bd') \mathbf{j} + (ad'+da'+bc'-cb') \mathbf{k }
\end{eqnarray*}
The conjugate of $ q = a \mathrm{I} + b \mathbf{i} + c \mathbf{j
}+ d \mathbf{k} $ is $ \overline{q} = a \mathrm{I} - b \mathbf{i}
- c \mathbf{j }- d \mathbf{k} $. A simple computation shows that
$\ds q \overline{q} = a^2 + b^2 + c^2 +d^2=\Delta$, $\ds q^{-1} =
\frac{1}{\Delta} \overline{q}$. The modulus of $q$ is $|q| =
\sqrt{\Delta}$. Since $ \mathbf{i}, \mathbf{j}, \mathbf{k} $ do
not commute, the multiplication of quaternions is noncommutative.

\section{$\mathbb{S}^3$ as a noncommutative Lie group}

Consider  $\mathbb{S}^3 = \{q \in \mathbf {Q}; |q| =1   \}$. For
any $q_1, q_2 \in \mathbb{S}^3$, $|q_1^{-1} q_2 | = \frac{\ds|
q_2|}{\ds |q_1|}=1$, {\it i.e.} $\mathbb{S}^3$ is a subgroup of
$\mathbf {Q}$. The law group on $\mathbb{S}^3$ is
\begin{eqnarray*}
&&(a,b,c,d) * (a', b', c', d')\\
& =&  (aa' - bb' - cc' - dd', ab'+ba' +cd' - dc',  ac'+ca'+db' -
bd', ad' + da' + bc' - cb'),
\end{eqnarray*}
with the identity element $(1,0,0,0)$ and the inverse
$(a,b,c,d)^{-1} = (a, -b, -c, -d)$. In this way $(\mathbb{S}^3,*)$
becomes a compact, noncommutative Lie group.\\

We shall work out the left invariant vector fields with respect to
the left translation $L_q : \mathbb{S}^3 \rightarrow
\mathbb{S}^3$, $L_q h = q *h$. These are vector fields on
$\mathbb{S}^3$ such that
$$ d L_q (X_h) = X_{q*h}, \qquad \forall q, h \in \mathbb{S}^3.$$
The set of all left invariant vector fields on $\mathbb{S}^3$
$$\mathcal{L}(\mathbb{S}^3)= \{ X; dL_q X = X, \forall q \in \mathbb{S}^3  \}  $$
is called the Lie algebra of $\mathbb{S}^3$.
$\mathcal{L}(\mathbb{S}^3)$ is a subalgebra of $\mathcal{L}(Q)$,
{\it i.e.} it is a subset closed under the Lie bracket $[\;,\;]$.
\\

In order to find a basis for $\mathcal{L}(\mathbb{S}^3)$, we shall
find first a basis for $\mathcal{L}(\mathbf{Q})$. Let $X$ be a
left invariant vector field on $\mathbf{Q}$. Then $X_q = (dL_q)
X_e$, where $e = (1,0,0,0)$ is the identity of $\mathbf{Q}$. If
write
$$ X_q = \sum_{i=1}^4 X_q^i \frac{\partial}{\partial x_i},$$
the components are
$$ X_q^i = (d L_q) X_e (x_i) = X_e (x_i \circ L_q),\quad i=1, \dots, 4,$$
where $x_i$ is the $i$-th coordinate.\\

Next we shall compute $x_i \circ L_q$, $i=1,2,3,4$. Let $q =
(a,b,c,d)$ and  $q' = (a',b', c', d')$. Then
\begin{equation}\label{eq:Lq}
(x_1\circ L_q) (q') = x_1 (L_q q') = aa' - bb' - cc' - dd'.
\end{equation}
Using
$$ a = x_1(q),\quad   a'= x_1 (q'), \quad b = x_2(q), \quad  b'= x_2 (q'),$$
$$ c = x_3(q),\quad  c'= x_3 (q'), \quad d= x_4(q), \quad  d'= x_4 (q'),$$
equation (\ref{eq:Lq}) becomes
$$(x_1 \circ L_q ) (q') = x_1(q) x_1(q') - x_2 (q) x_2(q') - x_3 (q) x_3 (q') - x_4(q) x_4(q'). $$
Dropping the argument $q'$, yields
$$x_1 \circ L_q  = x_1(q) x_1 - x_2 (q) x_2 - x_3 (q) x_3  - x_4(q) x_4. $$
By a  similar computation we find
\begin{eqnarray*}
 x_2 \circ L_q &=& x_1(q) x_2 + x_2(q) x_1 + x_3 (q) x_4 - x_4(q)
x_3,\\
 x_3 \circ L_q &=&  x_1(q) x_3 + x_3 (q) x_1 + x_4(q) x_2 - x_2 (q)
x_4,\\
x_4 \circ L_q &=& x_1(q) x_4 + x_4 (q) x_1 + x_2 (q) x_3 - x_3 (q)
x_2.
\end{eqnarray*} The first component becomes
\begin{eqnarray*}
X_q^1 &=& X_e(x_1 \circ L_q) \\
&=& x_1(q) X_e (x_1) - x_2(q) X_e(x_2) - x_3 (q) X_e (x_3) -
x_4(q) X_e(x_4)\\
&=&x_1(q) X_e^1 - x_2(q) X_e^2 - x_3 (q) X_e^3 - x_4(q) X_e^4.
\end{eqnarray*}
Dropping $q$ and denoting $\xi^i = X_e^i$, $i =1,2,3,4$, the
previous relation yields
$$X^1 = x_1 \xi^1 - x_2 \xi^2 - x_3 \xi^3 - x_4 \xi^4. $$
In a similar way we obtain \begin{eqnarray*}
 X^2 &=& x_1 \xi^2 + x_2
\xi^1 + x_3 \xi^4 - x_4 \xi^3,\\
 X^3 &=& x_1 \xi^3 + x_3 \xi^1 + x_4 \xi^2 - x_2 \xi^4,\\
 X^4 &=& x_1 \xi^4 + x_4 \xi^1 + x_2 \xi^3 - x_3 \xi^2.
\end{eqnarray*}
Left invariant vector fields $Z$ can be written as
\begin{eqnarray*}
Z = \sum_{i=1}^4 Z^i \frac{\partial }{\partial x_i}&=& \Big( x_1
\xi^1 - x_2 \xi^2 - x_3 \xi^3 - x_4 \xi^4 \Big) \frac{\partial
}{\partial x_1}\\
&&+\Big( x_1 \xi^2 + x_2 \xi^1 +x_3 \xi^4 - x_4 \xi^3 \Big)
\frac{\partial
}{\partial x_2}\\
&&+\Big( x_1 \xi^3 + x_3 \xi^1 +x_4 \xi^2 - x_2 \xi^4 \Big)
\frac{\partial
}{\partial x_3}\\
&&+\Big( x_1 \xi^4 + x_4 \xi^1 +x_2 \xi^3 - x_3 \xi^2 \Big)
\frac{\partial
}{\partial x_4}.\\
\end{eqnarray*}
Reordering, yields
\begin{eqnarray*}
Z &=& \xi^1 \Big[x_1 \frac{\partial}{\partial x_1}  +
x_2\frac{\partial}{\partial x_2} + x_3 \frac{\partial}{\partial
x_3} + x_4 \frac{\partial}{\partial x_4} \Big]\\
&&+ \xi^2 \Big[-x_2 \frac{\partial}{\partial x_1}  + x_1
\frac{\partial}{\partial x_2} + x_4 \frac{\partial}{\partial
x_3} - x_3 \frac{\partial}{\partial x_4} \Big]\\
&&+  \xi^3 \Big[-x_3 \frac{\partial}{\partial x_1}  - x_4
\frac{\partial}{\partial x_2} + x_1 \frac{\partial}{\partial
x_3} + x_2 \frac{\partial}{\partial x_4} \Big]\\
&&+  \xi^4 \Big[-x_4 \frac{\partial}{\partial x_1}  + x_3
\frac{\partial}{\partial x_2} -x_2 \frac{\partial}{\partial
x_3} + x_1 \frac{\partial}{\partial x_4} \Big]\\
&=& \xi^1 N - \xi^2 X - \xi^3 T - \xi^4 Y,
\end{eqnarray*}
where
$$ N = \sum_{i=1}^4 x_i \frac{\partial }{\partial x_i} $$
is the normal vector field to $\mathbb{S}^3$, and
\begin{eqnarray*}
X&=& x_2 \frac{\partial}{\partial x_1} - x_1
\frac{\partial}{\partial x_2} - x_4 \frac{\partial}{\partial x_3}
+ x_3 \frac{\partial}{\partial x_4},\\
Y&=& x_4 \frac{\partial}{\partial x_1} - x_3
\frac{\partial}{\partial x_2} + x_2 \frac{\partial}{\partial x_3}
- x_1 \frac{\partial}{\partial x_4},\\
T&=& x_3 \frac{\partial}{\partial x_1} + x_4
\frac{\partial}{\partial x_2} - x_1 \frac{\partial}{\partial x_3}
- x_2 \frac{\partial}{\partial x_4}.
\end{eqnarray*}
Using the $4\times 4$-matrix representations:
\begin{equation}
\label{eq:matrix}
I_1=\left[\array{rrrr}0 & 1 & 0 & 0
\\ -1 & 0 & 0 & 0
\\
0 & 0 & 0 & -1
\\ 0 & 0 & 1 & 0\endarray\right],\quad
I_3=\left[\array{rrrr}0 & 0 & 1 & 0
\\ 0 & 0 & 0 & 1
\\
-1 & 0 & 0 & 0
\\ 0 & -1 & 0 & 0\endarray\right],\quad
I_2=\left[\array{rrrr}0 & 0 & 0 & 1
\\ 0 & 0 & -1 & 0
\\
0 & 1 & 0 & 0
\\ -1 & 0 & 0 & 0\endarray\right],
\end{equation}
and the standard inner product $\langle\,,\,\rangle$ on $\mathbb{R}^4$, we also can write
\begin{equation*}N=\langle xU,\nabla_x\rangle,\quad
X=-\langle xI_1,\nabla_x\rangle,\quad Y=-\langle xI_3,\nabla_x\rangle,\quad
T=-\langle xI_2,\nabla_x\rangle,
\end{equation*}
where $U$ is the $4\times 4$ identity matrix.

It is easy to see that $\{ N, X, Y, T\}$ is a basis for
$\mathcal{L}(\mathbf{Q})$. In order to find a basis for
$\mathcal{L}(\mathbb{S}^3)$ we shall consider the vector fields
which are tangent to $\mathbb{S}^3$. Since a computation
provides
$$ \langle N, X \rangle =  \langle N, Y  \rangle = \langle N, T  \rangle = 0,$$
it follows that $X, Y, T$ are tangent to $\mathbb{S}^3$. Hence any
left invariant vector field on $(\mathbb{S}^3, *)$ is a linear
combination of $X, Y$ and $T$. Since the matrix of coefficients
$$ \begin{pmatrix}
 x_2 & -x_1 & -x_4 & \;\;\;x_3 \\
  x_4 & -x_3 &\;\;\; x_2 & -x_1 \\
  x_3 & \;\;\;x_4 & -x_1 & -x_2
\end{pmatrix}$$
has rank 3 at every point, it follows that $X, Y$ and $T$ are
linear independent and hence form a basis of $L(\mathbb{S}^3)$. It
worth nothing that the above rank property is not true without the
constraint $x_1^2 + x_2^2 + x_3^3 + x_4^2 =1$.\\

We also note that $\{X, Y, T, N\}$ form an orthonormal system with
respect to the usual inner product of $\mathbb{R}^4$. Then
$\mathcal{L}(\mathbf{Q}) = \mathcal{L}(\mathbb{S}^3)\bigoplus
\mathbb{R}N$.\\

\subsection{The horizontal distribution} Let $\mathcal{H} = span
\{ X, Y\}$ be the distribution generated by the vector fields $X$
and $Y$. Since $[Y, X] =  2 T \notin \mathcal{H} $ it follows that
$\mathcal{H}$ is not involutive. We can write
$\mathcal{L}(\mathbb{S}^3) = \mathcal{H}  \bigoplus \mathbb{R} T$.
The distribution $\mathcal{H}$ will be called the {\it horizontal
distribution}. Any curve on the sphere which has the velocity
vector contained in the distribution $\mathcal{H}$ will be called
a {\it horizontal curve}.

\smallskip
\begin{re}
Notice that this situation differs from the Heisenberg group since
the choice of the horizontal distribution is not unique. Because
of $[X,T]=-2Y$ and $[Y,T]=2X$ we could chose $\mathcal{H} = span
\{ X, T\}$ or $\mathcal{H} = span \{ Y, T\}$. The geometries
defining by different horizontal distributions are cyclically
symmetric, so we restrict out attention to the $\mathcal{H} = span
\{ X, Y\}$.
\end{re}

The following result deals with a characterization of horizontal
curves. It will make sense for later reasons to rename the
variables $y_1 = x_3$, $y_2 = x_4$.

\begin{prop}
\label{prop:hor} Let $\gamma (s) = (x_1(s), x_2(s), y_1(s),
y_2(s))$ be a curve on $\mathbb{S}^3$. The curve $\gamma$ is
horizontal if and only if
$$ \langle \dot x, y \rangle  = \langle x, \dot y \rangle.$$
\end{prop}
\begin{proof}
Since $\{ X, Y, T\}$ span the tangent space of the sphere
$\mathbb{S}^3$ we have
$$\dot \gamma = a X + b Y + c T. $$
Using that $\{X, Y, T  \}$ is an orthonormal system on
$\mathbb{S}^3$, we have \begin{eqnarray*}
 c &=& \langle  \dot \gamma
, T \rangle\\
&=& \dot x_1 y_1 + \dot x_2 y_2 -\dot y_1 x_1 - \dot y_2 x_2\\
&=& \langle \dot x, y \rangle  - \langle x, \dot y \rangle.
\end{eqnarray*}
The curve $\gamma$ is horizontal if and only if $c =0$ {\it i.e.}
$\langle
\dot x, y \rangle  = \langle x, \dot y \rangle$.\\
\end{proof}

Hence a horizontal curve $\gamma=(x_1, x_2, y_1, y_2): (0, \tau)
\rightarrow \mathbb{R}^3 $ has the velocity $ \dot \gamma = a X +
b Y $ where
\begin{eqnarray}
 a &=& \langle \dot \gamma, X \rangle= \dot x_1 x_2 - \dot x_2 x_1
 - \dot y_1 y_2 + \dot y_2 y_1 \label{eq:a222}\\
 b&=&\langle \dot \gamma, Y \rangle= \dot x_1 y_2 - \dot x_2 y_1 +
 \dot y_1 x_2 - \dot y_2 x_1.\label{eq:b222}
\end{eqnarray}
If we consider a metric on ${\cal H}$ such that $\{Y, X \}$ are
orthonormal vector fields, then the length of the curve $\gamma$ is
$$ \ell(\gamma) = \int_0^T \sqrt{a(s)^2 + b(s)^2} \, ds.$$
The horizontal curves minimizing length will be treated in a next
section.

\smallskip In the rest of this section we shall make a few
considerations regarding the horizontal distribution given as the
kernel of a one-form.

\smallskip
Consider the following one-form on $\mathbb{R}^4$
\begin{equation}
\omega = x_1 dy_1 - y_1 dx_1 + x_2 dy_2 - y_2 dx_2  = x dy - y dx.
\end{equation}
 One can easily check that
\begin{eqnarray*}
\omega(X) &=& -x_1 y_2 - y_1 x_2 + x_2 y_1 + x_1 y_2 =0\\
\omega(Y) &=& x_1 x_2 - y_1 y_2 - x_1 x_2 + y_1 y_2=0\\
\omega (T) &=& -x_1^2 -y_1^2 -x_2^2 -y_2^2 =-1 \not=0\\
\omega (N) &=& x_1 y_1 - y_1 x_1 + x_2 y_2 - y_2 x_2 =0,
\end{eqnarray*}
where
$$N = x_1 \partial_{x_1} + x_2 \partial_{x_2} + y_1 \partial_{y_1} +y_2 \partial_{y_2}$$
is the unit normal to the sphere $\mathbb{S}^3$. Hence $\ker
\omega = span\{ X, Y, N \}$ and the horizontal distribution can be
written as
$$\mathbb{S}^3 \ni p \rightarrow \mathcal{H}_p = \ker \omega \cap T_p \mathbb{S}^3. $$
Hence a vector $\nu= (v, w) = (v_1, v_2, w_1,w_2) \in \mathcal{H}$
if and only if
\begin{eqnarray*}
\langle x,v \rangle + \langle y, w \rangle  &=& 0 \;  (\Longleftrightarrow \nu \perp N)\\
\langle  v, y\rangle - \langle x, w\rangle &=&  0\;
(\Longleftrightarrow \omega (\nu)
=0)\\
|x|^2 + |y|^2 &=& 1 (  \Longleftrightarrow (x, y) \in
\mathbb{S}^3).
\end{eqnarray*}

\section{Connectivity theorem on $\mathbb{S}^3$}

The present section deals with the global connectivity property by
horizontal curves. Even if a result of Chow [8] states the
connectivity by piece-wise horizontal curves, here we shall prove
a globally smooth version of it. More precisely, the main result
of this section states that given two points on $\mathbb{S}^3$,
there is a horizontal smooth curve joining them.

\smallskip
Before proceeding to the proof, we need to introduce some
terminology and notations. This deals with the parametrization of
$\mathbb{S}^3$ in terms of the Euler angles, which due to
spherical symmetry are more suitable than the cartesian
coordinates.

\smallskip Consider $\varphi, \psi, \theta$ to be the Euler's
angles and let
$$\alpha = \frac{\varphi + \psi}{2}, \quad \beta = \frac{\varphi - \psi}{2}. $$
The sphere $\mathbb{S}^3$ can be parametrized as
\begin{eqnarray*}
x_1 &=& \cos \frac{\varphi + \psi}{2} \cos \frac{\theta}{2} = \cos
\alpha \cos \frac{\theta}{2},\\
x_2 &=& \sin \frac{\varphi + \psi}{2} \cos \frac{\theta}{2} = \sin
\alpha \cos \frac{\theta}{2},\\
y_1 &=& \cos \frac{\varphi - \psi}{2} \sin \frac{\theta}{2} = \cos
\beta \sin \frac{\theta}{2},\\
y_2 &=& \sin \frac{\varphi - \psi}{2} \sin \frac{\theta}{2} = \sin
\beta \sin \frac{\theta}{2},\\
\end{eqnarray*}
with $ \quad 0 \leq
 \theta \leq \pi,\quad -\pi \leq \alpha \leq \pi$. In the
 following we shall write the restriction of the one-form
 $$ \omega = x_1 dy_1 - y_1 dx_1 + x_2 dy_2 - y_2 dx_2$$
 to $\mathbb{S}^3$ using Euler's angles.
Since
\begin{eqnarray*}
dx_1 &=& -\sin \alpha \cos \frac{\theta}{2} d \alpha - \frac{1}{2}
\cos \alpha \sin \frac{\theta}{2} d\theta\\
dx_2 &=& \;\;\;\cos \alpha \cos \frac{\theta}{2} d \alpha -
\frac{1}{2}
\sin \alpha \sin \frac{\theta}{2} d\theta\\
dy_1 &=& -\sin \beta \sin \frac{\theta}{2} d \beta + \frac{1}{2}
\cos \beta \cos \frac{\theta}{2} d \theta\\
dy_2 &=& \;\;\; \cos \beta \sin\frac{\theta}{2} d \beta +
\frac{1}{2} \sin \beta \cos \frac{\theta}{2} d \theta,
\end{eqnarray*}
we obtain
\begin{eqnarray*}
\omega &=& x_1 dy_1 - y_1 dx_1 + x_2 dy_2 - y_2 dx_2\\
&=& \;\;\; \Big( \cos \beta \sin \alpha \sin \frac{\theta}{2} \cos
\frac{\theta}{2} - \sin \beta \cos \alpha \sin \frac{\theta}{2}
\cos \frac{\theta}{2} \Big) d \alpha \\
&&+ \Big( -\cos \alpha \cos \frac{\theta}{2} \sin \beta \sin
\frac{\theta}{2} + \sin \alpha \cos \frac{\theta}{2} \sin
\frac{\theta}{2} \cos \beta \Big) d \beta\\
&&+\frac{1}{2} \Big( \cos \alpha \cos \frac{\theta}{2} \cos \beta
\cos \frac{\theta}{2} + \cos \alpha \cos \beta \sin^2
\frac{\theta}{2} + \sin \alpha \sin \beta \cos^2
\frac{\theta}{2}+\sin \beta \sin \alpha \sin^2 \frac{\theta}{2}
\Big)d \theta\\
&=& \sin \frac{\theta}{2} \cos \frac{\theta}{2} \sin (\alpha -
\beta) d \alpha + \sin \frac{\theta}{2} \cos \frac{\theta}{2} \sin
(\alpha - \beta) d \beta + \frac{1}{2} \cos (\beta - \alpha) d
\theta\\
&=&\frac{1}{2} \sin \theta \sin \psi (d \alpha + d \beta)
+\frac{1}{2} \cos \psi d \theta\\
&=& \frac{1}{2} \Big( \sin \theta \sin \psi\,d \varphi + \cos
\psi\, d \theta \Big).
\end{eqnarray*}
Hence
$$ \omega_{|\mathbb{S}^3} = \frac{1}{2} \Big( \sin \theta \sin \psi\, d \varphi + \cos
\psi\, d \theta \Big).$$ The constraint $\sin \theta \sin \psi d
\varphi + \cos \psi d \theta =0$ is nonholonomic since
$$ 2 d \omega = \cos \theta \sin \psi \,d\theta \wedge d \varphi +\sin \theta \cos \psi\, d \psi \wedge d \varphi
 -\sin \psi\, d \psi \wedge d \theta \not=0.$$

We finally arrived  at the following  characterization of the
horizontal curves using  Euler's angles.

\begin{lemma}
Let $c(s) = (\varphi(s), \psi(s), \theta(s))$ be a curve on
$\mathbb{S}^3$. The curve $c$ is horizontal if and only if $\omega
(\dot c ) =0$, {\it i.e.}
\begin{equation}\label{eq:sinth1}
\sin \theta \sin \psi \,\dot \varphi + \cos \psi \,\dot \theta =
0.
\end{equation}
\end{lemma}

Before getting into the proof of the connectivity theorem, which
is the main result of this section, a couple of lemmas are needed.
The first one is a standard Calculus exercise.

\begin{lemma}\label{l:sph1}
Given the numbers $\alpha, \beta, \gamma \in \mathbb{R}$, there is
a smooth function \\$f :[0,1] \rightarrow \mathbb{R}$ such that
$$f(0) =0, \qquad f(1) = \alpha, \qquad f'(0) = \beta, \qquad f'(1)
= \gamma.$$
\end{lemma}
$$\epsfxsize=3.7in \epsfbox{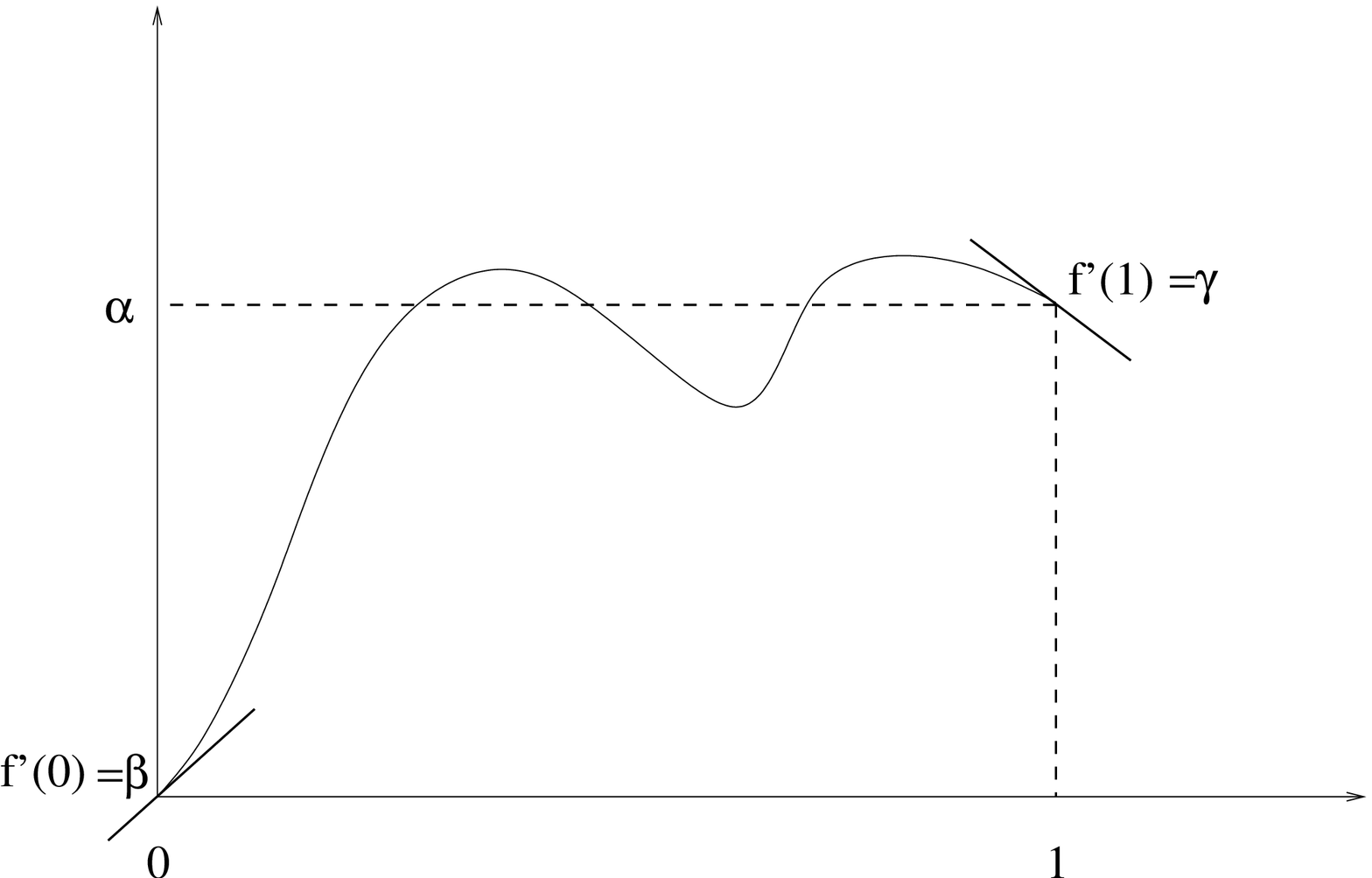}$$

\begin{lemma}\label{l:lem3}
Given $q_0, q_1 , I \in \mathbb{R}$, there is a function $q: [0,
1] \rightarrow \mathbb{R}$ such that
$$q(0) = q_0, \quad q(1) = q_1, \quad \int_0^1 q(u) \, du = I.  $$
\end{lemma}

\begin{proof}
Applying Lemma \ref{l:sph1}, there is a function $f :[0,1]
\rightarrow \mathbb{R}$ such that
$$ f(0) = 0, \quad f(1) = I, \quad f'(0) = q_0, \quad f'(1) = q_1.$$
Let $q(s) = f'(s)$. Then $f(s) = \int_0^s q(u) \, du$ and hence $I
= f(1) = \int_0^1 q(u) \, du$ and $q(0) = f'(0) = q_0$ and $q(1) =
f'(1) = q_1$.
\end{proof}

The next result is the theorem of connectivity announced at the
beginning of this section.
\begin{theorem}
Given two points $P, Q \in \mathbb{S}^3$, there is a smooth
horizontal curve joining $P$ and~$Q$.
\end{theorem}
\begin{proof}
Let $P(\varphi_0, \psi_0, \theta_0)$ and $Q(\varphi_1, \psi_1,
\theta_1)$ be the coordinates in the Eulerian angles. We shall
find a horizontal curve $c(s) = (\varphi(s), \psi(s), \theta(s))$
with $c(0) = P$ and $c(1) = Q$. This is equivalent with finding
functions $\varphi(s), \psi(s), \theta(s)$ such that
\begin{equation}\label{eq:shphor}
\sin \theta(s) \sin \psi(s) \dot \varphi(s) + \cos \psi(s) \dot
\theta(s) = 0,
\end{equation}
which satisfy the boundary conditions
$$\varphi(0) = \varphi_0, \qquad \psi(0) = \psi_0, \qquad \theta(0) = \theta_0, $$
$$\varphi(1) = \varphi_1, \qquad \psi(1) = \psi_1, \qquad \theta(1) = \theta_1.
$$
Assuming $\cos \psi(s) \not=0$, (\ref{eq:shphor}) can be written
as \begin{equation}\label{eq:shphor1}
 \sin \theta(s) \tan \psi(s) \dot \varphi(s) +
\dot \theta(s) =0.
\end{equation}
 Let
$$ \theta(s) = \arcsin p(s), \quad \psi(s) = \arctan q(s), \qquad \varphi(s) = \varphi_0 + s ( \varphi_1 - \varphi_0)$$
for some functions $p(s)$, $q(s)$ which will be determined later.
Let $k = \varphi_1 - \varphi_0$. Then (\ref{eq:shphor1}) yields
\begin{equation}
p(s) \, q(s) \, k + \frac{p'(s)}{\sqrt{1- p^2(s)}} =0.
\end{equation}
Separating and solving for $p(s)$ yields
$$ \frac{dp}{p\sqrt{1-p^2}} = - k q(s) \, ds$$
and integrating,
\begin{eqnarray*}
-\tanh^{-1}\frac{1}{\sqrt{1-p^2}} &=& -k \Bigg( \int_0^s q(u) \,
du + C_1 \Bigg) \Longleftrightarrow\\
\frac{1}{\sqrt{1-p^2(s)}} &=& \tanh  \Bigg[ k \Bigg( \int_0^s q(u)
\, du + C_1 \Bigg) \Bigg].
\end{eqnarray*}
The constant $C_1$ and the integral $\int_0^1 q(u) \, du$ can be
determined from the boundary conditions. We obtain
\begin{eqnarray*}
\frac{1}{\sqrt{1-p^2(0)}} = \tanh (k C_1)  &\Longrightarrow& C_1
=\frac{1}{k} \tanh^{-1} \Bigg(\frac{1}{\sqrt{1-p^2(0)}}\Bigg)\\
\frac{1}{\sqrt{1-p^2(1)}} = \tanh \Bigg[ k \Big(\int_0^1 q(u) \,
du + C_1  \Big) \Bigg] &\Longrightarrow&
\end{eqnarray*}
$$ \int_0^1 q(u) \, du = \frac{1}{k} \Bigg[ \tanh^{-1} \frac{1}{\sqrt{1-p^2(1)}} -
\tanh^{-1} \frac{1}{\sqrt{1-p^2(0)}}  \Bigg], $$ where $ p(0) =
\sin \theta(0)=\sin \theta_0 $ and $p(1) = \sin \theta_1$. Hence
$q(s)$ has to satisfy
\begin{equation}\label{eq:intb1}
 \int_0^1 q(u) \, du =
\frac{1}{k} \Bigg[ \tanh^{-1} \frac{1}{|\cos \theta_1|} -
\tanh^{-1} \frac{1}{|\cos \theta_0|}  \Bigg].
 \end{equation}
The boundary conditions also yield
\begin{equation}\label{eq:intb2}
 q(0) =\tan \psi_0, \qquad q(1) = \tan \psi_1.
 \end{equation}
Lemma \ref{l:lem3} provides the existence of a function
$q:[0,1]\rightarrow \mathbb{R}$ such that (\ref{eq:intb1}) and
(\ref{eq:intb2}) are satisfied.\\

The function $p(s)$ will be defined by the relation
$$ \frac{1}{\sqrt{1- p^2(s)}} = \tanh \Bigg[k \int_0^s q(u) \, du + \tanh^{-1} \frac{1}{|\cos \theta_0|}   \Bigg].$$
Hence the curve $ c(s) = (\varphi(s), \psi(s), \theta(s))=
(\varphi_0 + ks, \arctan q(s), \arcsin p(s))$ is the desired
horizontal curve joining the points $P$ and $Q$ on
$\mathbb{S}^3$.
\end{proof}

It is known that in the case of the Heisenberg group if $(x_0,0)$
and $(x_1,0)$ are two points in the $\{t=0\}$ plane, then the
segment joining the points is a horizontal curve joining the
points which lies in the $\{t=0\}$ plane. The following theorem is
an analog for the sphere $\mathbb{S}^3$. The plane will be
replaced by the sphere $\mathbb{S}^2$.

\begin{theorem}
Given $(\varphi_0, \varphi_1)$ and $(\theta_0, \theta_1)$, there
is a horizontal curve with $\psi =$ constant, which joins the
points with Euler coordinates $(\varphi_0, \psi, \theta_0)$ and
$(\varphi_1, \psi, \theta_1)$.
\end{theorem}
\begin{proof}
Let $(\varphi, \psi, \theta)$ be a horizontal curve with $\psi =$
constant. Then equation $(\ref{eq:sinth1}) $ can be written as
$$-\tan \psi \frac{d \varphi}{d \theta }= \frac{1}{\sin \theta}. $$
Integrating, yields
\begin{eqnarray}\label{eq:tanpsi}
-\tan \psi \int_{\theta_0}^{\theta} \frac{d \varphi}{d \theta } \,
d \theta &=& \int_{\theta_0}^{\theta} \frac{1}{\sin \theta} \, d
\theta  \Longleftrightarrow \nonumber \\
-\tan \psi \Big( \varphi(\theta) -\varphi(\theta_0) \Big) &=& \ln
\tan \frac{\theta}{2} - \ln \tan \frac{\theta_0}{2}.
\end{eqnarray}
For $\theta = \theta_1$ we obtain \begin{equation}\label{eq:psieq}
\psi = \arctan \Bigg({\ln \frac{\ds \tan \frac{\theta_1}{2}}{\ds
\tan \frac{\theta_0}{2}}}/(\varphi_0 - \varphi_1)  \Bigg).
 \end{equation}
Solving for $\varphi$, (\ref{eq:tanpsi}) yields
$$ \varphi(\theta) -\varphi_0 = - \ln \Big(\tan \frac{\theta}{2} / \tan \frac{\theta_0}{2} \Big) /\tan \psi, $$
with $\psi$ defined by (\ref{eq:psieq}). Hence the horizontal
curve joining $(\varphi_0, \psi, \theta_0)$ and $(\varphi_1, \psi,
\theta_1)$ is
$$(\varphi, \psi, \theta)=
\Big(  \varphi_0 -  \ln \Big(\tan \frac{\theta}{2} /
 \tan \frac{\theta_0}{2} \Big) /\tan \psi, \psi, \theta   \Big).  $$
\end{proof}

Let $S_{\psi_0} = \{ (\varphi, \psi, \theta)\in \mathbb{S}^3; \psi
= \psi_0 \} $ be the subset of $\mathbb{S}^3$ with the same Euler
angle $\psi$.  This corresponds to a 2-dimensional sphere.

\begin{corollary}
Given the points $(\varphi_0, \psi_0, \theta_0)$ and  $(\varphi_1,
\psi_1, \theta_1)$ on $\mathbb{S}^3$, with $\psi_0 = \psi_1$,
there is a horizontal curve in $S_{\psi_0}$  which joins the
points.
\end{corollary}

\section{The subRiemannian geodesics}

Given two points on the sphere $ \mathbb{S}^3$, we already know
that there are horizontal curves joining them. The horizontal
curve with the shortest length is called a geodesic. Standard
arguments show that an equivalent problem of finding the geodesics
 is to find the horizontal curves with minimum energy. This means
to minimize the action integral
$$ \int_0^T \frac{1}{2} [a^2(s) + b^2(s)] \, ds $$
subject to  the horizontality constraint $\langle \dot x, y
\rangle  = \langle x, \dot y \rangle$.  This will be minimizers
given by the solutions of the Euler-Lagrange system with
Lagrangian
$$ L(x, \dot x, y, \dot  y)  = \frac{1}{2} (a^2 + b^2)
 + \lambda(s) \big(x_1\dot y_1 - y_1 \dot x_1 + x_2 \dot y_2 - y_2 \dot x_2  \big), $$
where $a$ and $b$ depend on $x, \dot x, y, \dot y$, see formulas
(\ref{eq:a222})-(\ref{eq:b222}). The function $\lambda(s)$ is a
Lagrange multiplier function.
\medskip

\subsection{The case $\lambda(s) =0$.}

The curve which minimizes the energy in this case satisfies the
Euler-Lagrange system with the Lagrangian $\ds L(x, y, \dot x,
\dot y ) = \frac{1}{2} a^2 + \frac{1}{2} b^2$. We have
\begin{eqnarray*}
\frac{\partial L}{\partial \dot x_1} &=& ax_2 + b y_2, \qquad
\frac{\partial L}{\partial x_1} = -a \dot x_2 - b \dot y_2  \\
\frac{\partial L}{\partial \dot x_2} &=& -ax_1 - b y_1, \quad
\frac{\partial L}{\partial x_2} = a \dot x_1 + b \dot y_1\\
\frac{\partial L}{\partial \dot y_1} &=& -ay_2 + b x_2, \quad
\frac{\partial L}{\partial y_1} = a \dot y_2 - b \dot x_2\\
\frac{\partial L}{\partial \dot y_2} &=& ay_1 - b x_1, \quad
\frac{\partial L}{\partial y_2} = -a \dot y_1 +b \dot x_1.
\end{eqnarray*}

Then the Euler-Lagrange system
\begin{eqnarray}
&&\frac{d}{ds} \frac{\partial L}{\partial \dot x_1} =
\frac{\partial L}{\partial x_1}, \qquad \frac{d}{ds}
\frac{\partial L}{\partial \dot x_2} =
\frac{\partial L}{\partial x_2}, \label{eq:1stels} \\
&&\frac{d}{ds} \frac{\partial L}{\partial \dot y_1} =
\frac{\partial L}{\partial y_1}, \qquad \frac{d}{ds}
\frac{\partial L}{\partial \dot y_2} = \frac{\partial L}{\partial
y_2}, \label{eq:4stels}
\end{eqnarray}
becomes
\begin{eqnarray*}
&&\dot a x_2 + \dot b y_2 =  -2(a \dot x_2 + b \dot y_2),\qquad
 \dot a x_1 + \dot b y_1 = -2(a \dot x_1 + b \dot y_1),   \\
&&\dot a y_2 - \dot b x_2 = -2(a \dot y_2 - b \dot x_2), \qquad
 \dot a y_1 - \dot b x_1 = -2 (a \dot y_1 - b \dot x_1).
\end{eqnarray*}

Multiplying the first equation by $x_2$, the second by $x_1$, the
third by $y_2$ and the fourth by $y_1$, adding yields
\begin{eqnarray*}
&&\dot a (x_1^2 +x_2^2 + y_2^2+ y_1^2) + \dot b (y_2 x_2 + y_1 x_1 - x_2 y_2 - x_1
y_1)\\
&=& -2 [a(x_2 \dot x_2 + x_1 \dot x_1 + y_2 \dot y_2 + y_1 \dot
y_1) + b(x_2 \dot y_2 + x_1 \dot y_1 - y_2 \dot x_2 - y_1 \dot
x_1)] \\
\Longleftrightarrow\,\,\,\dot a &=& - 2b (y_2\dot  x_2 + y_1\dot
x_1 - x_2\dot y_2 - x_1
\dot y_1)\\
\Longleftrightarrow\,\,\,\dot a &=& 0,
\end{eqnarray*}
where we have used the horizontality condition (see Proposition
\ref {prop:hor})
$$ y_2\dot  x_2 + y_1\dot  x_1 - x_2\dot  y_2 - x_1
\dot y_1=0, $$
the constraint $$x_1^2 + x_2^2 + y_1^2 + y_2^2 =1$$
and its derivative
$$x_1 \dot x_1 + x_2 \dot x_2 + y_1 \dot y_1 +  y_2 \dot y_2 =0. $$

In a similar way, multiplying the first equation by $y_2$, the
second by $y_1$, the third by $-x_2$, and the fourth by $-x_1$,
adding, yields $\dot b =0$. We arrived at the following result,
which provides four constants of motion (three angles and the
energy).

\begin{theorem} Given two points $P, Q \in \mathbb{S}^3$, let
 $\gamma$ be a length minimizing curve among all the horizontal
curves which join $P$ and $Q$. Then
\begin{description}
  \item[1)] $|\dot \gamma(s) |$ is constant along the curve.
  \item[2)] The angles between the velocity $\dot \gamma$ and the vector
  fields $X$, $Y$ and $T$ are constant along the curve.
\end{description}
\end{theorem}

\begin{proof}
\begin{description}
  \item[1)] From $\dot a = \dot b =0$ it follows that $a(s) =
  $constant and $b(s) =
  $constant and hence
  $\ds |\dot \gamma(s)| = \sqrt{a(s)^2 + b(s)^2}=$ constant.
  \item[2)] Using $\bf 1)$, we have
  $$ \cos  \widehat{\dot \gamma, X }  = \frac{\langle \dot \gamma, X\rangle }{|\dot \gamma|\, |X|}
  =\frac{a(s)}{\sqrt{a(s)^2 + b(s)^2}}=constant $$
   $$ \cos  \widehat{\dot \gamma, Y }  = \frac{\langle \dot \gamma, Y \rangle }{|\dot \gamma|\, |Y|}
  =\frac{b(s)}{\sqrt{a(s)^2 + b(s)^2}}=constant $$
   $$ \cos  \widehat{\dot \gamma, T }  = \frac{\langle \dot \gamma, T \rangle }{|\dot \gamma|\, |T|}
  =0 \Longrightarrow  \widehat{\dot \gamma, T }  = \pi/2. $$
\end{description}
\end{proof}

$$\epsfxsize=3in \epsfbox{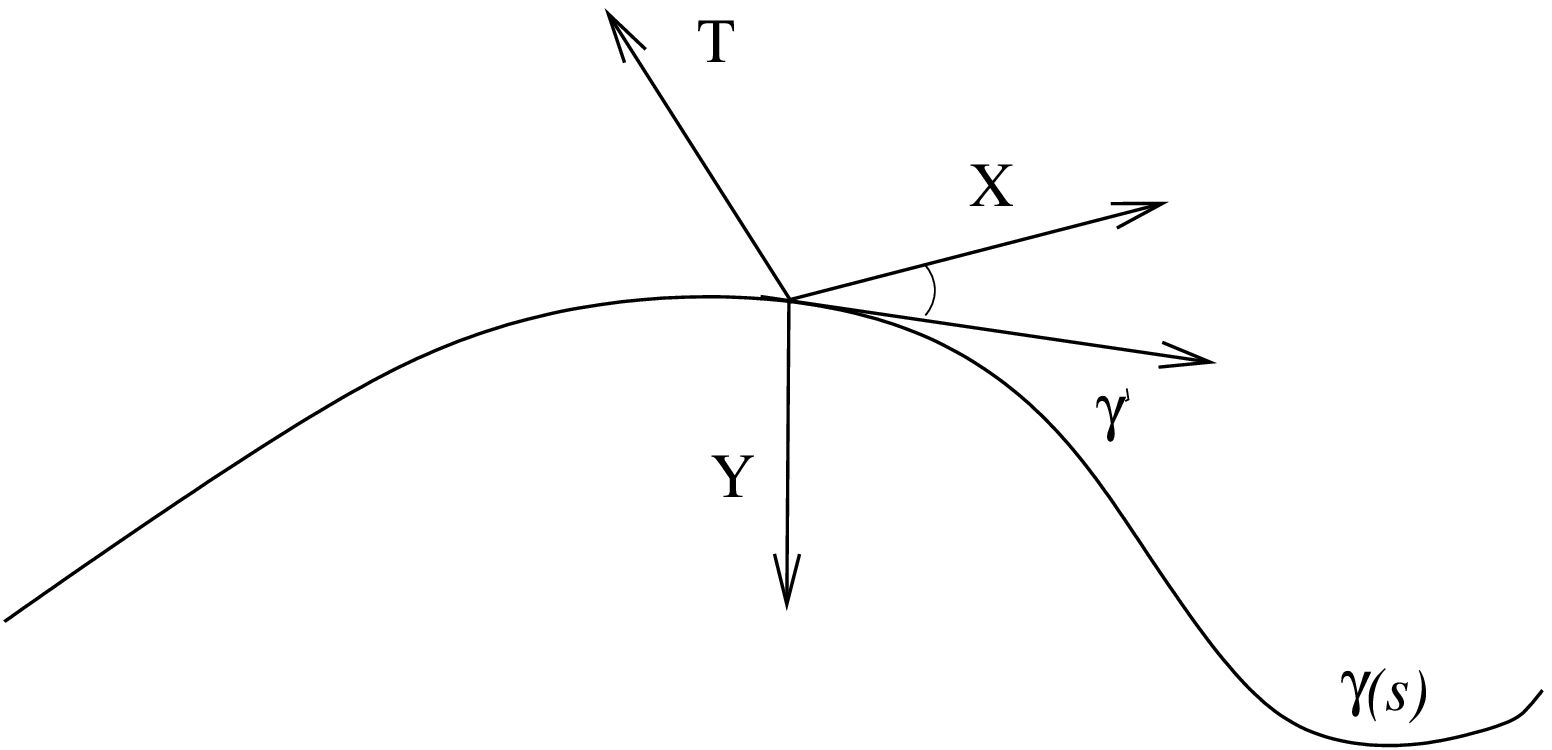}\quad
\epsfxsize=1.8in \epsfbox{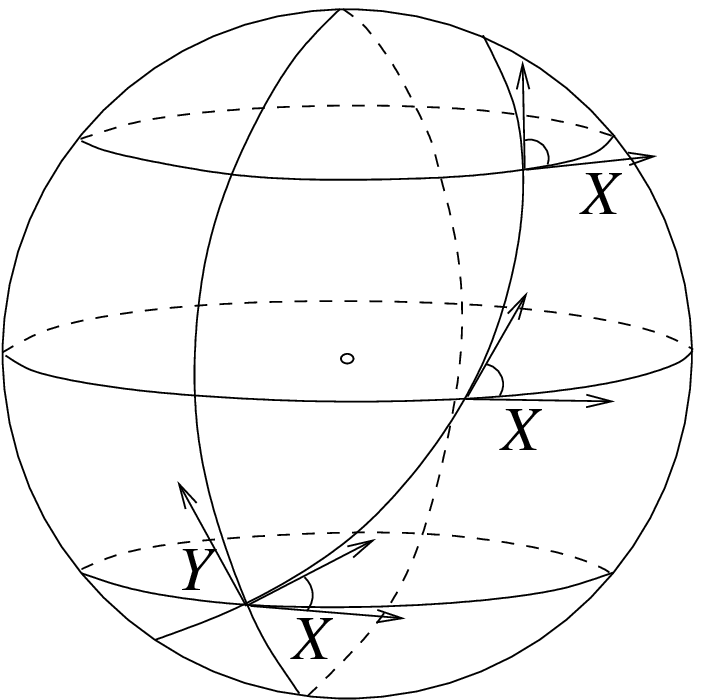}$$

 \centerline{\small {\it
Figure 1 : The angle between the velocity $\dot \gamma(s)$ and the
vectors $X$ and $Y$ is constant. }}

\medskip
\begin{re}
This case when the  Lagrange multiplier vanishes yields in the
case of the Heisenberg group lines parallel to the $x$-space.
These are particular cases of subRiemannian geodesics.
\end{re}
\smallskip

\subsection{The case $\lambda(s) \not=0$.}

The Lagrangian in this case is
$$L = \frac{1}{2} a^2 + \frac{1}{2} b^2 + \lambda (\dot x_1 y_1 + \dot x_2 y_2 - \dot y_1 x_1 - \dot y_2 x_2). $$
The Euler-Lagrange system (\ref{eq:1stels}) - (\ref{eq:4stels})
becomes
\begin{eqnarray*}
\dot a x_2 + \dot b y_2 &=& - 2 (a\dot x_2 + b \dot y_2 + \lambda
\dot y_1 ) - \dot \lambda y_1\\
\dot a x_1 + \dot b y_1 &=& - 2 (a\dot x_1 + b \dot y_1 - \lambda
\dot y_2 ) + \dot \lambda y_2\\
\dot a y_2 - \dot b x_2 &=& - 2 (a\dot y_2 - b \dot x_2 + \lambda \dot x_1) - \dot \lambda x_1\\
\dot a y_1 - \dot b x_1 &=& -2(a\dot y_1 - b \dot x_1 - \lambda
\dot x_2) + \dot \lambda x_2.
\end{eqnarray*}
Multiplying the first equation by $x_2$, the second by $x_1$, the
third by $y_2$ and the forth by $y_1$, adding yields
\begin{eqnarray*}
&&\dot a (\underbrace{x_1^2 + x_2^2 + y_1^2 + y_2^2}_{=1}) + \dot
b(\underbrace{y_2 x_2 + y_1 x_1 - x_2 y_2 - x_1 y_1}_{=0})\\
&=& -2 \Big( a (\underbrace{x_2 \dot x_2 + x_1 \dot x_1 + y_2 \dot
y_2 + y_1 \dot y_1}_{=0}) + b (\underbrace{x_2 \dot y_2 + x_1 \dot
y_1 - y_2 \dot x_2 - y_1 \dot x_1}_{=0}) \\
&&+\lambda (\underbrace{x_2 \dot y_1 - x_1 \dot y_2 + \dot x_1 y_2
- \dot x_2 y_1 }_{=b}) \Big) + \dot \lambda (\underbrace{- x_2 y_1
+ x_1 y_2 - x_1 y_2 + x_2 y_1}_{=0} ) \Longleftrightarrow\\
\dot a &=& - 2 \lambda b.
\end{eqnarray*}

In a similar way, multiplying the first equation by $y_2$, the
second by $y_1$, the third by $-x_2$ and the forth by $-x_1$,
adding we obtain

\begin{eqnarray*}
&& \dot a (\underbrace{x_2 y_2 + x_1 y_1 - x_2 y_2 - x_1
y_1}_{=0}) +
\dot b (\underbrace{y_2^2 + y_1^2 + x_2^2 + x_1^2}_{=1}) \\
&=& -2 \Big(  a (\underbrace{\dot x_2 y_2 + \dot x_1 y_1 - \dot
y_2 x_2 - \dot y_1 x_1}_{=0} ) + b (y_2 \dot y_2 + y_1\dot y_1 +
x_2 \dot x_2 + x_1 \dot x_1 )\\
&&+ \lambda ( \underbrace{\dot y_1 y_2 - \dot y_2 y_1 - \dot x_1
x_2 + \dot x_2 x_1}_{=-a}) + \dot \lambda (\underbrace{- y_1 y_2 +
y_2 y_1 - x_1 x_2 - x_2 x_1}_{=0})\Longleftrightarrow\\
\dot b&=& 2\lambda a.
\end{eqnarray*}
Hence $a$ and $b$ satisfy the following system

\begin{eqnarray}
\dot a &=& - 2 \lambda b \label{eq:aba1}\\
\dot b&=& 2\lambda a.\label{eq:aba2}
\end{eqnarray}
Multiplying the first equation by $a$ and the second by $b $ and
add, yields $a\dot a + b \dot b =0$, {\it i.e.}, $ a^2 + b^2 =
r^2$  constant along the geodesics. This means that the energy is
preserved and the velocity of the geodesics have constant length.
Let $\theta(s)$ be a function such that
$$ a(s) = r \cos \theta(s), \quad b(s) = r \sin \theta(s).$$
Substituting in the equations (\ref{eq:aba1}), (\ref{eq:aba2}), we
obtain
\begin{eqnarray*}
\sin \theta(s) \, \dot \theta(s) &=& 2 \lambda(s) \sin \theta(s)\\
\cos \theta(s) \, \dot \theta(s) &=& 2 \lambda(s) \cos \theta(s),
\end{eqnarray*}
which after dividing the first equation by $\sin \theta$ and
the second by $\cos \theta$ yields
$$ \dot \theta(s) = 2 \lambda(s). $$
The solution is $\theta(s) = 2 \Lambda(s) + \theta_0$, where $\ds
\Lambda (s) = \int_0^s \lambda(u) \, du $. Hence
\begin{eqnarray}
a(s) &=& r \cos \big( 2 \Lambda(s) + \theta_0 \big), \label{eq:ala1}\\
b(s) &=& r \sin \big( 2 \Lambda(s) + \theta_0
\big).\label{eq:ala2}
\end{eqnarray}
 Since $\Lambda(0) =0$,
then $a(0) = r \cos \theta_0$ and $b(0) = r \sin \theta_0$, which
provides $\theta_0 = \tan^{-1} (b(0)/a(0))$.

\medskip
Let $\gamma(s)$ be a subRiemannian geodesic.
 Since $\{X, Y\}$ are orthonormal,
 $$a =\langle\dot \gamma, X\rangle = \underbrace{|\dot \gamma|}_{=r} \cdot \underbrace{|X|}_{=1}
  \cos (\widehat{\dot \gamma, X})=   r \cos (\widehat{\dot \gamma, X})$$
and
 $$b =\langle\dot \gamma, Y\rangle = \underbrace{|\dot \gamma|}_{=r} \cdot \underbrace{|Y|}_{=1}
  \cos (\widehat{\dot \gamma, Y})=   r \cos (\widehat{\dot \gamma, Y})=
  r \sin \big(\frac{\pi}{2} - \widehat{\dot \gamma, Y}  \big).$$
Comparing with (\ref{eq:ala1}), (\ref{eq:ala2}) we arrive at the
following result.

\begin{proposition}\label{p:dotsR1} The angles under which the subRiemannian geodesics
intersect the integral curves of the vector fields $X$ and $Y$ are
given by the formulas
$$\widehat{ \dot \gamma(s), X_{\gamma(s)}}= 2 \Lambda(s) + \theta_0,\quad
 \widehat{ \dot
\gamma(s), Y_{\gamma(s)}} = \frac{\pi}{2}- \widehat{ \dot
\gamma(s), X_{\gamma(s)}}.  $$
\end{proposition}

The main goal now is to find $\Lambda(s)$. In order to do this we
shall construct an equivalent variational problem. We need the
following result, which writes the energy in a more friendly way.

\begin{prop}\label{prop5.4}
If $\big(x_1(s), x_2(s), y_1(s), y_2(s)\big)$ is a subRiemannian
geodesic then we have  $$\dot x_1^2(s) + \dot x_2^2 (s) + \dot
y_1^2(s) + \dot y_2^2(s)=a^2(s) + b^2(s). $$
\end{prop}
\begin{proof}
Using the definitions of $a$ and $b$, the horizontality condition,
and the holonomic constraint, we can write the following system

\begin{eqnarray*}
\dot x_1 x_2 + y_1 \dot y_2 - \dot x_2 x_1 - \dot y_1 y_2 &=&
a(s)\\
\dot x_1 y_2 + x_2 \dot y_1 - \dot x_2 y_1 - \dot y_2 x_1 &=&
b(s)\\
\dot x_1 y_1 + \dot x_2 y_2 - x_2 \dot y_2 - x_1 \dot y_1 &=& 0\\
\dot x_1 x_1 + \dot x_2 x_2 + \dot y_1 y_1 + \dot y_2 y_2 &=& 0.
\end{eqnarray*}
This can be written in a matrix way as
$$
\underbrace{\begin{pmatrix}
  x_{2} & -x_{1} & -y_{2} & y_{1} \\
  y_{2} & -y_{1} & x_{2} & -x_{1} \\
  y_{1} & y_{2} & -x_{1} & -x_{2} \\
  x_{1} & x_{2} & y_{1} & y_{2}
\end{pmatrix}}_{=M}\begin{pmatrix}
  \dot x_{1} \\
  \dot x_{2} \\
  \dot y_{1} \\
  \dot y_{2}
\end{pmatrix} = \begin{pmatrix}
  a \\
  b \\
  0 \\
  0
\end{pmatrix}.
$$
Since  $\det M =1 $, $M^{-1} = M^t$,
 $M$ is an orthogonal matrix,  and hence
preserves the  Euclidean length of vectors. It follows that
$$ \dot x_1^2 + \dot x_2^2 + \dot y_1^2 + \dot
y_2^2 = a^2 + b^2$$
 along any subRiemannian geodesic, which proves the statement.
\end{proof}

\begin{re}
The above result holds because the Euclidean metric on
$\mathbb{R}^4$ restricted on $span\{X, Y\}$ is the subRiemannain
metric.
\end{re}

The next result will be useful later in the sequel, when we shall
find the Lagrange multiplier.

\begin{prop}\label{p:glam6}
Let $\gamma(s)$ be a subRiemannain geodesic. Then $\langle \ddot
\gamma(s), T_{\gamma(s)} \rangle =0 $, for any $s$, {\i.e.}, the
component of the acceleration along the missing direction $T=[X,
Y]$ vanishes.
\end{prop}
\begin{proof} Let $\gamma(s) = \big(x_1(s), x_2(s), y_1(s),
y_2(s)\big)$. Differentiating in the horizontality condition
$$ \dot x_1 y_1 + \dot x_2 y_2  - \dot y_1 x_1 - \dot y_2 x_2 = 0$$
yields $\ddot x_1 y_1 + \ddot x_2 y_2 - \ddot y_1 x_1 - \ddot y_2
x_2 = 0 $, which can be written as
$$\langle (\ddot x_1, \ddot x_2, \ddot y_1, \ddot y_2), (y_1, y_2, -x_1, -x_2) \rangle =0,$$
or $\langle   \ddot \gamma, T \rangle =0 $, where
$$ T = y_1 \partial_{x_1} + y_2 \partial _{x_2} - x_1 \partial_{y_1} - x_2 \partial_{y_2}. $$
\end{proof}

We shall consider another variational problem, where we replace
the energy $\ds  \frac{1}{2}(a^2 + b^2)$ by $\ds \frac{1}{2} (\dot
x_1^2 + \dot x_2^2 + \dot y_1^2 + \dot y_2^2 )$ and consider the
new Lagrangian
$$L^* (x, \dot x, y, \dot y) =  \frac{1}{2} (\dot
x_1^2 + \dot x_2^2 + \dot y_1^2 + \dot y_2^2 ) +\lambda(s) (\dot
x_1 y_1 + \dot x_2 y_2 - \dot y_1 x_1 - \dot y_2 x_2 ).$$ The
Euler-Lagrange equations provided by the Lagrangian $L^*$ describe
the same subRiemannian geodesics as the equations associated with
the initial Lagrangian $L$. The Euler-Lagrange equations for $L^*$
are
\begin{eqnarray*}
\ddot x_1 &=& 2 \lambda \dot y_1 + \dot \lambda y_1\\
\ddot x_2 &=& 2 \lambda \dot y_2 + \dot \lambda y_2\\
\ddot y_1 &=& - 2 \lambda \dot x_1 - \dot \lambda x_1 \\
\ddot y_2 &=& - 2 \lambda \dot x_2 - \dot \lambda x_2.
\end{eqnarray*}

Multiplying the first equation by $y_1$, the second by $y_2$, the
third by $-x_1$, and the fourth by $-x_2$, adding yields
\begin{eqnarray*}
y_1 \ddot x_1 + y_2 \ddot x_2 - x_1 \ddot y_1 - x_2 \ddot y_2 &=&
2\lambda (\underbrace{y_1 \dot y_1 + y_2 \dot y_2 + x_1 \dot x_1 +
x_2 \dot x_2}_{=0}) \\&&+ \dot \lambda (\underbrace{y_1^2 + y_2^2
+ x_1^2 + x_2^2}_{=1}),
\end{eqnarray*}
which can be written as
\begin{eqnarray*}
\dot \lambda(s) &=&y_1 \ddot x_1 + y_2 \ddot x_2 - x_1 \ddot y_1 -
x_2 \ddot y_2 \\
&=& \langle (\ddot x_1, \ddot x_2, \ddot y_1, \ddot y_2), (y_1,
y_2, -x_1, -x_2) \rangle\\
&=& \langle \ddot \gamma, T \rangle.
\end{eqnarray*}

Using Proposition \ref{p:glam6} we obtain the following result.

\begin{prop}
The Lagrange multiplier $\lambda(s)$ is constant along any
subRiemannian geodesics.
\end{prop}

If let $\lambda = c/2$, with $c$ constant, then Proposition
\ref{p:dotsR1} yields
$$ \widehat{\dot \gamma(s), X_{\gamma(s)}} =  cs + \theta_0.$$
This leads to the following characterization of the subRiemannian
geodesics.

\begin{prop}
The  unit speed curve  $\gamma(s)$ is a subRiemannian geodesic on
$(\mathbb{S}^3, X, Y)$ if and only if the angle between its
velocity and the direction of the vector field $X$ increases
linearly in~$s$.
\end{prop}
\smallskip

\subsection{The group $SU(2)$} The special unitary group
 $$ SU(2) = \{x_1+ x_2 \mathbf{i} + y_1 \mathbf{j} +  y_2 \mathbf{k};\ \ 
 x_1^2 + x_2^2 + y_1^2 + y_2^2 = 1 \}
 $$
is also a Lie group, which can be identified with the sphere
$\mathbb{S}^3$ by the isomorphism
 $$ \varphi : SU(2) \rightarrow \mathbb{S}^3, \qquad \varphi(x_1I+ x_2\mathbf{i}+
 y_1\mathbf{j}
 + y_2\mathbf{k}) = (x_1, x_2, y_1, y_2). $$

Since
\begin{eqnarray}
\label{eq:5.19}
X_q &=& x_2 \partial_{x_1} -x_1 \partial_{x_2} -y_2 \partial_{y_1}
+ y_1 \partial_{y_2}\nonumber \\
&=& (x_2,  -x_1,   -y_2, y_1 )\cdot \nabla_q
= (x_2-x_1\mathbf{i} -y_2\mathbf{j}+ y_1\mathbf{k})\cdot \nabla_q \nonumber \\
&=& -(x_1I+ x_2\mathbf{i} + y_1\mathbf{j}+y_2\mathbf{k})\mathbf{i}\cdot \nabla_q
= -(x_1, x_2, y_1, y_2)\mathbf{i}\cdot \nabla_q \nonumber \\
& =& - \langle(q\mathbf{i}), \nabla_q\rangle,
\end{eqnarray}

\begin{eqnarray}\label{eq:5.20}
Y_q &=& y_2 \partial_{x_1} -y_1 \partial_{x_2} +x_2 \partial_{y_1}
-x_1 \partial_{y_2}\nonumber \\
&=&   (y_2,  -y_1, x_2,  -x_1  )\cdot \nabla_q
=  (y_2I  -y_1\mathbf{i}+x_2\mathbf{j} -x_1\mathbf{k})\cdot \nabla_q\nonumber \\
&=&-(x_1I  + x_2\mathbf{i}+ y_1\mathbf{j}+ y_2\mathbf{k})\mathbf{k}\cdot \nabla_q
= -(x_1, x_2, y_1, y_2)\mathbf{k}\cdot \nabla_q\nonumber \\
 & = & - \langle(q\mathbf{k}),\nabla_q\rangle,
\end{eqnarray}
\begin{eqnarray}
\label{eq:5.21}
T_q &=& y_1 \partial_{x_1} + y_2 \partial_{x_2}
- x_1 \partial_{y_1}  -x_2 \partial_{y_2}\nonumber \\
&=&  (y_1,  y_2, -x_1,  -x_2  )\cdot \nabla_q
=  (y_1I + y_2\mathbf{i} -x_1\mathbf{j}-x_2\mathbf{k})\cdot \nabla_q\nonumber \\
&=&-(x_1I +x_2\mathbf{i}+ y_1\mathbf{j}+ y_2\mathbf{k})\mathbf{j}\cdot \nabla_q
= -(x_1, x_2, y_1, y_2)\mathbf{j}\cdot \nabla_q\nonumber \\
 & = & -\langle(q\mathbf{j}),\nabla_q\rangle,
\end{eqnarray}
where $q=(x_1,x_2,y_1,y_2)$,
$\nabla_q=(\partial_{x_1},\partial_{x_2},\partial_{y_1},\partial_{y_2})$
and $\mathbf{i}$, $\mathbf{j}$, $\mathbf{k}$ are defined in~\eqref{eq:matrixun}. It follows
that the left invariant vector fields $X, Y, T$, which span the Lie
algebra of $\mathbb{S}^3$ correspond on $SU(2)$ to $-\mathbf{i},\, -\mathbf{k},\,
-\mathbf{j}$, respectively. \smallskip

\noindent {\bf Remark.} An alternative approach to calculate
geodesics is to use the {\it Hamiltonian method}. This method is
very complicated in the case of studying geodesics on ${\mathbb
S}^3$, but some information we can obtain from the superficial
analysis. Using the
notations~\eqref{eq:5.19},~\eqref{eq:5.20},~\eqref{eq:5.21}, the
vector fields $X,Y$ and $T$ can be written in the form
$$X=\langle -qI_1,\nabla_q\rangle, \quad Y=\langle -qI_3,\nabla_q\rangle,
\quad T=\langle -qI_2,\nabla_q\rangle
$$
where $\langle\cdot,\cdot\rangle$ is the usual scalar product in
$\mathbb R^4$. Define the Hamilton function
$$H=\frac{1}{2}(X^2+Y^2)=\frac{1}{2}\Big(\langle qI_1,\xi\rangle^2
+\langle qI_3,\xi\rangle^2\Big).$$ Then the Hamilton system is the
following
\begin{equation}\label{hs}
\begin{split}
&\dot q=\frac{\partial H}{\partial \xi}\ \ \Rightarrow \ \ \dot
q=\langle qI_1,\xi\rangle\cdot (qI_1)+\langle qI_3,\xi\rangle\cdot (qI_3)\\
&\dot \xi=-\frac{\partial H}{\partial q}\ \ \Rightarrow \ \
\dot\xi=\langle qI_1,\xi\rangle\cdot (\xi I_1)+\langle
qI_3,\xi\rangle\cdot (\xi I_3).
\end{split}\end{equation}
A geodesic is the projection of the solution of Hamilton's system to
the $q$-space. Since $\langle qI_1,q\rangle=\langle qI_2,
q\rangle=\langle qI_3,q\rangle=0$, multiplying the first equation
of~\eqref{hs} by $q$ we get
$$\langle\dot q,q\rangle=0\ \ \Rightarrow \ \ |q|^2=const.$$
We conclude that the solution of the Hamiltonian system belongs to a
sphere. Taking the constant equals $1$ we get geodesics belonging to
$\mathbb S^3$. Multiplying the first equation of~\eqref{hs} by
$qI_2$, we get
\begin{equation}\label{hc}
\langle\dot q,qI_2\rangle=0,
\end{equation}
by the role of multiplication between $I_1$, $I_2$, and $I_3$. The
reader easily recognize the horizontality condition $\langle\dot
x,y\rangle=\langle x,\dot y\rangle$ in~\eqref{hc}. It means that the
solution of the Hamiltonian system is a horizontal curve.
Multiplying the first equation of~\eqref{hs} by $qI_1$ or by $qI_3$,
we get
$$\langle \xi,qI_1\rangle=\langle\dot q,qI_1 \mathbf{i}\rangle,
\qquad \langle \xi,qI_3\rangle=\langle\dot q,qI_3\rangle.$$ We see
that the Hamilton function can also be written in the form
$$H=\frac{1}{2}\Big(\langle qI_1,\xi\rangle^2+\langle qI_3,\xi\rangle^2\Big)
=\frac{1}{2}\Big(\langle qI_1,\dot q\rangle^2+\langle qI_3,\dot
q\rangle^2\Big).
$$
If we multiply the first equation of~\eqref{hs}
by $\dot q$ then we get
\begin{equation}\label{en}
|\dot q|^2=\langle qI_1,\xi\rangle^2+\langle qI_3,\xi\rangle^2
=\langle qI_1,\dot q\rangle^2+\langle qI_3,\dot q\rangle^2=2H.
\end{equation}
Thus the Hamiltonian function give the kinetic energy $H=\frac{|\dot
q|^2}{2}$ and it is a constant along the geodesics. Notice that
$\langle qI_1,\dot q\rangle=a$ and $\langle qI_3,\dot q\rangle=b$
by~\eqref{eq:a222} and~\eqref{eq:b222}. We conclude $|\dot
q|^2=a^2+b^2$ from~\eqref{en} that corresponds to
Proposition~\ref{prop5.4}. This is only the beginning of a long term
project. There are many problems remain open. For example, given any
two points on ${\mathbb S}^3$, how many geodesics connecting them?
Is there any abnormal minimizer in this case? What is the action
function? What is the volume element which is the solution of a
transport equation? We will answer these questions in a forthcoming
paper.

\noindent {Department of Mathematics, Eastern Michigan University,
Ypsilanti, MI, 48197, USA}

\noindent {e-mail}: ocalin@emunix.emich.edu

\noindent {Department of Mathematics, Georgetown University,
Washington D.C., 20057, USA}

\noindent {e-mail}: chang@math.georgetown.edu

\noindent {Department of Mathematics, University of Bergen,
Johannes Brunsgate 12, Bergen 5008, Norway}

\noindent {e-mail}: irina.markina@uib.no

\end{document}